\documentclass[english]{article}
\usepackage[T1]{fontenc}
\usepackage[latin9]{inputenc}
\usepackage{amsmath}
\usepackage{graphicx}
\usepackage{amssymb}
\usepackage{esint}



\usepackage{babel}

\begin{document}

\title{Introduction to generalised C\'{e}saro convergence II}

\author{Richard Stone}
\maketitle
\begin{abstract}
\begin{eqnarray*}
\begin{array}{cc}
Beauty \; is \; truth, \; truth \; beauty \; - \; that \; is \; all \\
Ye \; know \; on \; earth, \; and \; all \; ye \; need \; to \; know. \\
\end{array}
\end{eqnarray*}

In this second of three introductory papers, we extend the notion of generalised C\'{e}saro summation/convergence to the more natural setting of what we call \textit{remainder} C\'{e}saro summation/convergence. This greatly expands the range of problems susceptible to C\'{e}saro methods and introduces the \textit{geometric} location of summands as a critical consideration. We also show that geometric generalised C\'{e}saro convergence is invariant under dilation and scaling. We present a number of calculations illustrating the utility of these developments. In particular we introduce a new, more natural definition of the classical Gamma function using remainder C\'{e}saro summation/products, and show that many its key properties - both basic and advanced - fall out directly and intuitively from this C\'{e}saro definition and its geometric and dilation-invariance properties. We also consider other examples and show how C\'{e}saro methodology explains the common structure of many well-known functional equations.
\end{abstract}

\section{Introduction}

This paper extends the theory of generalised C\'{e}saro convergence introduced in [1].
 
In section 2 we introduce a core new notion - remainder C\'{e}saro summation/convergence, where the partial-sum function analysed via a C\'{e}saro scheme is now on a contour relative to chosen point $z_{0}\in\mathbb{C}$, rather than always being on $[0,\infty)$. We give both a formal and a working definition of this expanded notion, and discuss how the geometric location of the summands now becomes critical - so much so that we often refer henceforth to  \textit{geometric} generalised C\'{e}saro summation/convergence.

The move to remainder C\'{e}saro summation/convergence greatly increases the range of problems we can tackle and techniques we can apply. To begin with, in subsection 2.1 we see that it allows us to make sense of finite sums, $\sum_{j=1}^{k}f(j)$, where $k$ no longer need be a discrete positive integer, but can instead be an arbitrary complex parameter. As an illustration we show that the Hurewicz zeta function, $\zeta_{H}$, is defined most naturally in this context. From this we further show that a well-known differential equation relating the polynomials in $k$ given by $\sum_{j=1}^{k}j^{n}$, $n\in\mathbb{Z}_{\geq{0}}$ (which are closely related to the Bernoulli polynomials) follows naturally within this remainder C\'{e}saro approach. In the same vein, we show how to reformulate the intriguing integral identity from [1, section 4.3], namely that
\begin{equation}
\int_{-1}^{0}\,\sum_{j=1}^{k}j^{n}\,dk\,=\,\zeta(-n)
\label{zeta_Integral_Identity_1}\end{equation}
in a way that applies for all $\rho\in\mathbb{C}$, not just $\rho=n\in\mathbb{Z}_{\geq{0}}$, albeit that we will have to defer the proof of this more general identity until after we have discussed C\'{e}saro dilation-invariance in section 3.\footnote{Overall, the discussion in section 2.1 effectively answers the question posed in [1] about why discrete variables are often treatable as continuous, subject to techniques of calculus.}

In section 2.2 we extend remainder C\'{e}saro summation to consider not just summation from $z_{0}$ "in the positive direction", but also "in the negative direction" and in both directions ('bi-directional"); and we consider corresponding notions of remainder C\'{e}saro products. All of these notions effectively define operators on suitable spaces of complex functions. Applying the remainder C\'{e}saro product to the elementary case of $f(z)=z$ gives a function, ${\tilde{\Gamma}}(z+1)$, which formally emulates the factorial function for $z\in\mathbb{Z}_{\geq{0}}$ and clearly should give $\Gamma(z+1)$ in a more direct and natural way than any of the definitions traditionally found in textbooks and the literature. Understanding ${\tilde{\Gamma}}(z+1)$ provides an almost-perfect test-lab\footnote{Labradors are excellent dogs} for understanding the new ideas and techniques of geometric remainder C\'{e}saro methodology, and so we devote the final subsection of section 2 to analysing it in detail. 

Section 2.3 contains three distinct calculations for this C\'{e}saro emulation of the Gamma function, each showing a key feature of the new geometric remainder C\'{e}saro methodology. The first shows that ${\tilde{\Gamma}}(z+1)$ is in fact well-defined as a function on all of $\mathbb{C\setminus}\mathbb{Z}_{\leq{0}}$ - and in the course of the proof the criticality of the geometric placement of summands, as well as of the geometric treatment of divergences, within the remainder C\'{e}saro framework becomes clear.

The second then calculates ${\tilde{\Gamma}}\,^{\prime}\left(0\right)$ and the higher derivatives, ${\tilde{\Gamma}}^{(n)}\left(0\right)$, $n\in\mathbb{Z}_{\geq{2}}$, allowing us to deduce the Taylor series for ${\tilde{\Gamma}}$ near $0$ and from it the equivalence of ${\tilde{\Gamma}}(z+1)$ and ${\Gamma}(z+1)$ on all of $\mathbb{C\setminus}\mathbb{Z}_{\leq{0}}$. The calculation of ${\tilde{\Gamma}}^{\prime}\left(0\right)$ is particularly interesting. It again demonstrates that it is essential to pay careful heed to the geometric location of summands when performing generalised remainder C\'{e}saro calculations, and by so doing achieves what, in other Physics-related contexts, would be understood as a renormalisation calculation.

Finally, having established that $\tilde{\Gamma}$ is in fact just a more natural, C\'{e}saro definition of $\Gamma$, we derive its famous functional equation that
\begin{equation}
\Gamma(z)\,\Gamma(1-z)\,=\,\frac{\pi}{\sin\left(\pi z\right)}
\label{GammaFunctionalEquation}\end{equation}
using this C\'{e}saro definition. The key to the proof is that bi-directional summation automatically sends any suitable input function to a periodic function with period $1$, which allows us to invoke Fourier theory. Much of the basic structure of this proof is universal, and thus equally applicable to other functions defined via remainder C\'{e}saro summation, not just $\Gamma$. As such this proof in turn sheds light on why the functional equations of many other famous functions likewise relate their function values at $z$ to their values at $1-z$, and we consider some such examples.

In section 3 we then turn to the behaviour of geometric C\'{e}saro convergence under two important transformation groups - dilation and scaling.

First, in section 3.1 we show that geometric C\'{e}saro convergence is dilation invariant. We give two proofs. The first is direct. For the second, we recall how useful it was to introduce the inverse operator, $P_{D}^{-1}$, when working within the discrete C\'{e}saro framework in [1, section 3.4] and so calculate the corresponding inverse continuous C\'{e}saro operator $P^{-1}$. We then deduce dilation-invariance based on the operator commutation relation between the generator of dilations and $P^{-1}$, which is trivial to calculate. This approach has a dual value - it shows how closely related $P^{-1}$ (and hence $P$) is to the group of dilations and its generator (which is interesting in and of itself), and it indicates the approach we will follow later in proving that geometric C\'{e}saro convergence is also scaling-invariant.

Before turning to this scaling-invariance, however, we use section 3.2 to demonstrate that C\'{e}saro dilation-invariance is not merely an abstract curiosity - it leads immediately to clean, intuitive proofs of interesting, non-trivial results. Specifically, returning to the Gamma function, we consider one of its more advanced properties, namely its multiplication formulae. The proof of these in textbooks is often rather obscure, giving little insight into \textit{why} they are true and how they were originally deduced (see e.g. [2]). We show that in the geometric remainder C\'{e}saro formulation they become completely natural, and their proof becomes an almost trivial consequence of C\'{e}saro dilation-invariance. As before, the essential elements of the proof are applicable equally well to other functions defined via remainder C\'{e}saro summation - not just to $\Gamma$ - so that analogous formulae\footnote{Although in general these are duplication \textit{summation} formulae rather than \textit{multiplication} formulae as for $\Gamma$; the $\Gamma$ multiplication formulae are really summation formulae for $\ln\left(\Gamma\right)$.} apply to many other such well-known functions. As an example, we derive the analogous result for the Hurewicz zeta function, $\zeta_{H}$. 

In section 3.3, we prove that C\'{e}saro convergence is scaling-invariant.

Finally, in section 4 we provide appendices where we fill in the working for some claims made in section 2. Specifically, we provide the derivation of the functional equation for $\zeta$ using bi-directional summation and Fourier theory; and we show how the remainder C\'{e}saro perspective leads immediately to efficient code for calculation of $\Gamma(z+1)$ for arbitrary $z\in\mathbb{C}$.

\section{Remainder C\'{e}saro summation/convergence}

In [1] the notion of (continuous) C\'{e}saro convergence was recast in terms of a C\'{e}saro operator $P[f](x):=\frac{1}{x}\intop_{0}^{x}f(t)\,\textrm{d}t$ (on functions $f:\mathbb{R}_{>0}\rightarrow\mathbb{C}$) and extended by consideration of the eigenvalues and eigenfunctions of $P$, so that a classically divergent function $f$ has generalised C\'{e}saro limit $L$ if $q(P)[f](x)\rightarrow L$ classically as $x\rightarrow\infty$ for some regular polynomial $q(P)$ (regularity being equivalent to having $q(1)=1$). This allows us to annihilate (i.e. assign generalised C\'{e}saro limit $0$ to) divergent eigenfunctions and generalised eigenfunctions of $P$ with eigenvalues $\lambda\neq1$ - namely powers $x^{\rho}$, and power-log functions $x^{\rho}\left(ln x\right)^{m}$ for $\rho\in\mathbb{C}\setminus\left\{0\right\}$, $m\in\mathbb{Z}_{\geq{0}}$ - before trying to average the residual component of $f$ sufficiently often (by applying $P^{n}$ for sufficiently large $n$) to make it convergent to limit $L$. Since these eigenfunctions $x^{\rho}$ arise naturally in the Euler-McLaurin asymptotic expansion for the p-sum function\footnote{i.e. partial sum function} of the Riemann zeta series, $\sum_{n=1}^{\infty}n^{-s}$, we saw in [1, section 3] that generalised C\'{e}saro summation can be used to derive the analytic continuation of $\zeta$ to $Re(s)\leq1$ and understand its singularities.

In applying this geometric C\'{e}saro approach to the treatment of $\zeta$, the p-sum function (of the real variable $x:=k+\alpha$) is obtained by adding in each term $n^{-s}$ at the point $z=n$. We can view this as a special case of what we shall call \textit{remainder} C\'{e}saro summation in which contributions occur at points located relative to an initial point $z_{0}\in\mathbb{C}$.

Specifically, given a function $f:\mathbb{C\rightarrow\mathbb{C}}$ we wish to define its strict remainder sum at $z_{0}\in\mathbb{C}$ by
\begin{equation}
R_{+}[f](z_{0}):=\sum_{n=1}^{\infty}f(z_{0}+n)
\label{eq:Remsum}\end{equation}
where we understand the sum on the RHS in a generalised C\'{e}saro sense along the horizontal contour $\gamma:\, t\mapsto z_{0}+t$. We do this by adapting the working definition in theorem 2 in [1] into the following working definition of C\'{e}saro convergence along a contour:\\
\\
\textbf{Definition 1:} \textit{Suppose} $\gamma:[0,\infty)\rightarrow\mathbb{C}:\, t\mapsto\gamma(t)$ \textit{is a contour parametrised by arc length} $t$, \textit{starting at} $z_{0}$ \textit{with} $\gamma(t)\rightarrow\infty$ \textit{as} $t\rightarrow\infty$, \textit{and suppose} $f$ \textit{is a function on} $\gamma$ \textit{which can be written as}
\begin{equation*}
f(t)=\sum_{j=1}^{n}a_{j}(\gamma(t))^{\rho_{j}}(\ln(\gamma(t)))^{m_{j}}+R(t)
\end{equation*}
\textit{for some finite collection of constants} $a_{j}\in\mathbb{C}$, $\rho_{j}\in\mathbb{C\setminus}\{0\}$ \textit{and} $m_{j}\in\mathbb{Z_{\geq\textrm{0}}}$ \textit{and some remainder function} $R(t)$. \textit{If there exists} $n\in\mathbb{Z_{\geq\textrm{0}}}$ \textit{such that} $P^{n}[R](t)\rightarrow L$ \textit{as} $t\rightarrow\infty$, \textit{then we say that} $f$ \textit{has generalised C\'{e}saro limit} $L$ \textit{along} $\gamma$, \textit{where} $P$ \textit{is now the averaging operator} (\textit{in} $t$) \textit{along} $\gamma$ \textit{defined by}
\begin{equation*}
P[h](t)=\frac{1}{t}\int_{0}^{t}h(u)\,\textrm{d}u
\end{equation*}
\\
\textbf{Notes on Definition 1:} \textbf{(i)} In the rest of this paper the only contours we shall need to consider are either horizontal rays parallel to the real axis ($\gamma:\, t\mapsto z_{0}\pm t$) or vertical rays parallel to the imaginary axis ($\gamma:\, t\mapsto z_{0}\pm it$), which simplifies interpretation. For these cases it is not hard to see uniqueness of limits in definition 1.\\
\\
\textbf{(ii) [Criticality of geometry]} It is critical in this definition, and will be crucial throughout this paper, that the functions $z^{\rho_{j}}(\ln z)^{m_{j}}$ that we throw away (i.e. assign generalised C\'{e}saro limit $0$ for $\rho_{j}\neq0$) are functions of the \textit{geometric} variable $z=\gamma(t)$, not just of the arc-length parameter $t$. We continue to think of these as C\'{e}saro eigenfunctions and generalised eigenfunctions in this setting even though this association is now less straightforward given the definition of $P$ in terms of $t$ rather than $z$.\footnote{In fact we can continue to think of $P$ as given by $P[f](z)=\frac{1}{z}\int_{0}^{z}\,f(w)\,dw$ so that these are still eigenfunctions and generalised eigenfunctions, but we omit further discussion here.} In the case where $\rho_{j}\notin\mathbb{Z}$ the need for this distinction disappears, at least when $\gamma$ is a simple ray as per (i), since in this case a straightforward Taylor expansion allows us to express each $z^{\rho_{j}}(\ln z)^{m_{j}}$ as a linear combination of functions of the form $t^{\tilde{\rho}_{i}}(\ln t)^{\tilde{m}_{i}}$ with all $\rho_{i}\neq0$, and vice-versa. However, when $\rho_{j}$ is an integer this distinction becomes pivotal in order to avoid pure powers of $\ln$, which correspond to C\'{e}saro eigenfunctions with eigenvalue $1$ and thus have no generalised C\'{e}saro limit (see [1], section 2).\\
\\
\textbf{(iii)} In relation to this last remark, recall the key point that C\'{e}saro definition of generalised convergence is fundamentally intended as a tool for constructive analytic continuation - in some complex parameter $s\in\mathbb{C}$ which drives the values of the $\rho_{j}$'s in the situation being analysed. As such, the above distinction could be rephrased as saying that at a point $s$ where all $\rho_{j}\notin\mathbb{Z}$ we do not need to be careful to distinguish between the parameter $t$ and the geometric variable $z=\gamma(t)$; but when $s$ is such that there is a $\rho_{j}\in\mathbb{Z_{\geq\textrm{0}}}$ we will need to distinguish between $z$ and $t$ as per the discussion in (ii) in order to obtain the correct analytic continuation across this value of $s$.\\
\\
\textbf{(iv)} When $\gamma$ is a horizontal ray and $z_{0}=0$, definition 1 devolves simply to the case considered in [1] for the analytic continuation of $\zeta(s)$; so definition 1 represents a natural extension of the generalised C\'{e}saro definition in [1], allowing us now also to consider horizontal rays starting at arbitrary $z_{0}\in\mathbb{C}$, or vertical rays or indeed more general contours $\gamma$. Note, however, that while the motivating case of remainder summation above leads to a partial sum function with evenly spaced jumps, in definition 1 the function $f$ need not have any jumps or even arise from a summation process; or even if it does, it may arise in a way where the spacing between summands is not even.\\
\\
\textbf{(v)} We will return to some of these observations at various points during computations later in this paper, and their significance may then become clearer. For now, however, we simply conclude by summarising definition 1 into the following working recipe for calculation of generalised C\'{e}saro limits along contours:\\

\textbf{(a)} First remove linear combinations of C\'{e}saro eigenfunctions and generalised eigenfunctions, but critically doing so geometrically as eigenfunctions in $z=\gamma(t)$, rather than simply in $t$, and then

\textbf{(b)} Apply a suitable power of the C\'{e}saro averaging operator along the contour (i.e. averaging in the contour arc-length $t$).\\
\\
In the case of remainder summation this means we form the contour $\gamma:\, t\mapsto z_{0}+t$, form the p-sum function $s_{f}(z_{0};z)=s_{f}(z_{0};z_{0}+k+\alpha):=\sum_{n\leq k}\, f(z_{0}+n)$ and use definition 1 to obtain
$R_{+}[f](z_{0})$ as $\underset{z\rightarrow\infty}{Clim}\, s_{f}(z_{0},z)$. Here we have used the notation $\underset{z\rightarrow\infty}{Clim}$ rather than $\underset{t\rightarrow\infty}{Clim}$ or $\underset{k\rightarrow\infty}{Clim}$ to emphasise that it is eigenfunctions and generalised eigenfunctions in $z=z_{0}+k+\alpha$, rather than in $t$ or $k$, which we remove from $s_{f}(z_{0},z)$ prior to averaging as discussed above.

\subsection{The Hurewicz zeta function; Example calculations; Finite sums and discrete variables revisited}

Let us illustrate remainder C\'{e}saro summation with some example calculations, starting with the case of $f(z)=z^{-s}$. \\
\\
\textbf{Hurewicz zeta function:} In this case we get a function which is already well-known, the Hurewicz zeta function
\begin{equation}
\zeta_{H}\left(z_{0};s\right):=R_{+}\left[\tilde{z}^{-s}\right]\left(z_{0}\right) \quad .
\label{Hurewicz_zeta}\end{equation}
This is now a function of two complex variables, $z_{0}$ and $s$; and again is classically convergent for $Re(s)>1$ (at least for $z_{0}\notin\mathbb{Z}_{<0}$) and coincides with $\zeta(s)$ when $z_{0}=0$ (i.e. $\zeta(s)=\zeta_{H}\left(0;s\right)=R_{+}\left[\tilde{z}^{-s}\right]\left(0\right)$ for all $s\in\mathbb{C\setminus}\{1\}$)\footnote{The tilde notation, $\tilde{z}^{-s}$, in definition \ref{Hurewicz_zeta} means of course, $R_{+}[f](z_{0})$ where $f(z)=z^{-s}$, with the remainder summation in the $z$-variable for given, fixed $s$; we shall use this notational shorthand extensively and even sometimes write $R_{+}^{\left(\tilde{z}\right)}$ for clarity.}.

As with the treatment of $\zeta(s)$ in [1], the analytic continuation of $\zeta_{H}\left(z_{0};s\right)$ to $Re(s)\leq 1$ for arbitrary fixed $z_{0}$ then proceeds strip-wise, first to $0<Re(s)\leq 1$, then to $-1<Re(s)\leq 0$ etc, but now using definition 1 in concert with the Euler-McLaurin sum formula, and using the binomial expansion to re-write terms of the form $\left(z_{0}+k\right)^{\rho}$ in the p-sum function in descending powers ($z^{\rho}$, $z^{\rho-1}$ etc) of the geometric variable $z=z_{0}+k+\alpha$. As example calculations, let us consider the two cases, $s=0$ and $s=-1$, and now calculate $\zeta_{H}\left(z_{0};0\right)$ and $\zeta_{H}\left(z_{0};-1\right)$ for arbitrary $z_{0}$, just like we calculated $\zeta(0)$ and $\zeta(-1)$ in [1].\\
\\
\textbf{(i)} When $s=0$ we have that $s_{f}\left(z_{0},0;z\right)=s_{f}\left(z_{0},0;z_{0}+k+\alpha\right)=\sum_{j=1}^{k}1=k=\left(z_{0}+k+\alpha\right)-z_{0}-\alpha=z-z_{0}-\alpha$. Thus
\begin{equation}
\zeta_{H}\left(z_{0};0\right)=R_{+}\left[\tilde{z}^{0}\right]\left(z_{0}\right)=-z_{0}-\frac{1}{2}
\label{Hurewicz_zeta_0}\end{equation}
on noting that $z=z_{0}+k+\alpha$ is an eigenfunction of $P$ with eigenvalue $\frac{1}{2}$ and that we still have $P\left[\tilde{\alpha}\right]\left(z\right)=\frac{1}{k+\alpha}\{\sum_{j=0}^{k-1}\int_{0}^{1}\,\tilde{\alpha}\,d\tilde{\alpha}+\int_{0}^{\alpha}\,\tilde{\alpha}\,d\tilde{\alpha}\}=\frac{1}{k+\alpha}\{\frac{1}{2}k+\frac{1}{2}\alpha^{2}\}\rightarrow\frac{1}{2}$ as $k\rightarrow\infty$.\\
\\
\textbf{(ii)} In the same way, for $s=-1$ we have $s_{f}\left(z_{0},-1;z\right)=\sum_{j=1}^{k}(z_{0}+j)=z_{0}k+\frac{1}{2}k^{2}+\frac{1}{2}k=\frac{1}{2}\left(z_{0}+k+\alpha\right)^{2}-\frac{1}{2}z_{0}^{2}-\frac{1}{2}\alpha^{2}-\left(z_{0}+k\right)\alpha+\frac{1}{2}k=\frac{1}{2}z^2-\frac{1}{2}z_{0}^{2}-\frac{1}{2}z_{0}-z_{0}\left(\alpha-\frac{1}{2}\right)-k\left(\alpha-\frac{1}{2}\right)-\frac{1}{2}\alpha^{2}$. Now, we saw in (i) that $\underset{z\rightarrow\infty}{Clim}\,\alpha=\frac{1}{2}$, so that $\underset{z\rightarrow\infty}{Clim}\,z_{0}\left(\alpha-\frac{1}{2}\right)=0$; and $P\left[\tilde{k}\left(\tilde{\alpha}-\frac{1}{2}\right)\right]\left(z\right)=\frac{1}{k+\alpha}\{\sum_{j=0}^{k-1}j\int_{0}^{1}\,\left(\tilde{\alpha}-\frac{1}{2}\right)\,d\tilde{\alpha}+k\int_{0}^{\alpha}\,\left(\tilde{\alpha}-\frac{1}{2}\right)\,d\tilde{\alpha}\}=\frac{k}{k+\alpha}\{\frac{1}{2}\alpha^{2}-\frac{1}{2}\alpha\}=\frac{1}{2}\alpha^{2}-\frac{1}{2}\alpha+o(1)$. Thus since, as above, $\underset{z\rightarrow\infty}{Clim}\,\alpha=\frac{1}{2}$ and $\underset{z\rightarrow\infty}{Clim}\,\alpha^{2}=\frac{1}{3}$, it follows that we have
\begin{eqnarray}
\zeta_{H}\left(z_{0};-1\right)=R_{+}\left[\tilde{z}^{1}\right]\left(z_{0}\right) & = & \underset{z\rightarrow\infty}{Clim}\,\{\frac{1}{2}z^2-\frac{1}{2}z_{0}^{2}-\frac{1}{2}z_{0}-\alpha^{2}+\frac{1}{2}\alpha\}\\
 & = & -\frac{1}{2}z_{0}^{2}-\frac{1}{2}z_{0}-\frac{1}{12}.
\label{Hurewicz_zeta_-1}\end{eqnarray}
\\
\textbf{(iii)} In both cases (i) and (ii) we could of course guess the form of $\zeta_{H}\left(z_{0};0\right)$ and $\zeta_{H}\left(z_{0};-1\right)$ heuristically by first considering $z_{0}=k\in\mathbb{Z}_{>0}$. In this case we have $R_{+}\left[\tilde{z}^{-s}\right]\left(k\right)=R_{+}\left[\tilde{z}^{-s}\right]\left(0\right)-\sum_{j=1}^{k}j^{-s}$ and since $R_{+}\left[\tilde{z}^{0}\right]\left(0\right)=\zeta(0)=-\frac{1}{2}$ and $\sum_{j=1}^{k}1=k$ we have that $\zeta_{H}\left(k;0\right)=-k-\frac{1}{2}$; while for $s=-1$ we have $R_{+}\left[\tilde{z}^{1}\right]\left(0\right)=\zeta(-1)=-\frac{1}{12}$ and $\sum_{j=1}^{k}j=\frac{1}{2}k^{2}+\frac{1}{2}k$ so that $\zeta_{H}\left(k;-1\right)=-\frac{1}{2}k^{2}-\frac{1}{2}k-\frac{1}{12}$. The formulae for arbitrary $z_{0}$ in equations \ref{Hurewicz_zeta_0} and \ref{Hurewicz_zeta_-1} then follow immediately by simply substituting $z_{0}$ in place of $k$ in the obvious way.\\
\\
\textbf{Extending finite sums from discrete to continuous variables:} In fact, remainder C\'{e}saro summation using definition 1 is precisely what makes sense of this heuristic working and explains why, in many cases, we can treat the discrete variable, $k$, in a finite sum, $\sum_{j=1}^{k}$, as if it were a continuous variable.

Specifically, any finite sum $\sum_{j=1}^{k}f(j)$, viewed as a function of the discrete variable $k\in\mathbb{Z}_{>0}$, extends naturally to a function of $z_{0}\in\mathbb{C}$ by defining $\sum_{j=1}^{z_{0}}f(j)$ as
\begin{equation}
\sum_{j=1}^{z_{0}}f(j):=R_{+}\left[f\right]\left(0\right)-R_{+}\left[f\right]\left(z_{0}\right) \quad .
\label{Finite_sum_z_0}\end{equation}
We may then differentiate or integrate this function with respect to $z_{0}$, just as though the original variable $k$ had been a complex variable.\\
\\
\textbf{Examples: (i)} For example, extending the calculations for $s=0$ and $s=-1$ above, suppose we consider the family of polynomials in $k$, $b_{n}(k)$, defined by $b_{n}(k):=\sum_{j=1}^{k}j^{n-1}$ for $n\in\mathbb{Z}_{>0}$. Then we find that these are related by formal differentiation of the discrete summation variable, $k$, namely 
\begin{equation}
\frac{d}{dk}(b_{n}(k))=(n-1)\cdot\left\{b_{n-1}(k)-\zeta(2-n)\right\}
\label{b_n_deriv}\end{equation}
from which they may be readily recursively generated and shown to be connected to the Bernoulli polynomials $B_{n}(x)$. Why does this work? It is because, defining $b_{n}(z)$ for arbitrary $z\in\mathbb{C}$ as above by
\begin{equation}
b_{n}(z):=R_{+}\left[\tilde{z}^{n-1}\right]\left(0\right)-R_{+}\left[\tilde{z}^{n-1}\right]\left(z\right)=\zeta(-n+1)-\zeta_{H}\left(z;-n+1\right)
\label{b_n_z_0}\end{equation}
we may rigorously differentiate w.r.t $z$. Since $\frac{d}{dz}\zeta_{H}(z,s)=-s\zeta_{H}(z,s+1)$ (after noting that summation and differentiation may be commuted in a careful C\'{e}saro analysis) we then immediately obtain the claimed relationship on setting $z=k$.\\
\\
\textbf{(ii)} In the same way, in [1, section 4.3] we proved the intriguing integral relationship given in equation \ref{zeta_Integral_Identity_1}, which involves integration w.r.t a discrete summation variable, $k$, namely that
\begin{equation*}
\int_{-1}^{0}\,\sum_{j=1}^{k}j^{n}\,dk\,=\,\zeta(-n)
\end{equation*}
for any $n\in\mathbb{Z}_{\geq{0}}$. How to make sense of this? Well, understanding $\sum_{j=1}^{z}j^{n}$ as $R_{+}\left[\tilde{z}^{n}\right]\left(0\right)-R_{+}\left[\tilde{z}^{n}\right]\left(z\right)=\zeta(-n)-\zeta_{H}\left(z;-n\right)$, this is equivalent to saying that
\begin{equation}
\int_{-1}^{0}\,\zeta_{H}\left(z;-n\right)\,dz\,=\,0 \quad \textrm{for all} \quad n\in\mathbb{Z}_{\geq{0}}
\label{zeta_Integral_Identity_2}\end{equation}
and this now makes perfect sense as an integral in a continuous variable $z$. In fact, this re-formulation of the result allows us to loosen the requirement that $n\in\mathbb{Z}_{\geq{0}}$ and conjecture that in general,
\begin{equation}
\int_{-1}^{0}\,\zeta_{H}\left(z;s\right)\,dz\,=\,0 \quad \textrm{for any} \quad s\in\mathbb{C\setminus}\{1\}.
\label{zeta_Integral_Identity_3}\end{equation}

In [1] the proof of equation \ref{zeta_Integral_Identity_1} used the C\'{e}saro asymptotic relationship that $\underset{z\rightarrow\infty}{Clim}\,k^{n}=\frac{(-1)^{n}}{n+1}$ for $n\in\mathbb{Z}_{\geq{0}}$\footnote{Here $k$ is measured relative to $z_{0}=0$, as was the case throughout [1]}, and also relied on having $n\in\mathbb{Z}_{\geq{0}}$ to ensure that $\sum_{j=1}^{k}j^{n}$ is a polynomial in $k$ with only a finite number of integer powers of $k$ to which this result can then be applied. It turns out that conjecture \ref{zeta_Integral_Identity_3} is in fact true for arbitrary $s\in\mathbb{C\setminus}\{1\}$. 

Its proof for non-integer $s$, however, relies on the dilation-invariance of generalised geometric C\'{e}saro convergence. As such we content ourselves here with having explained how to make sense of the apparent integration w.r.t a discrete variable in equation \ref{zeta_Integral_Identity_1} and we defer the proof of equation \ref{zeta_Integral_Identity_3} until section 3.

\subsection{Remainder C\'{e}saro operators; Remainder C\'{e}saro products; A C\'{e}saro definition of the Gamma function}

We can extend remainder C\'{e}saro summation by considering summation not just to "the right", but also to "the left" and in both directions, since definition 1 certainly allows us to consider also contours which are horizontal rays in the negative real direction. We define $R_{+,0}[f](z):=\sum_{n=0}^{\infty}f(z+n)$, $R_{-}[f](z):=\sum_{n=1}^{\infty}f(z-n)$, $R_{-,0}[f](z):=\sum_{n=0}^{\infty}f(z-n)$, $R_{+,0,-}[f](z):=\sum_{n=-\infty}^{\infty}f(z+n)$ and so forth. In the last case with a bi-directional sum, the definition is via two independent C\'{e}saro sums, for $R_{+,0}$ and $R_{-}$ respectively.

Note that $R_{+}$, $R_{-}$, $R_{+,0,-}$ etc are now properly understood as operators on suitable spaces of functions of a complex variable. We do not explore this in any detail here\footnote{Indeed neither the author nor his cats have undertaken any such systematic exploration of the operator and spectral properties of $R_{+}$, $R_{-}$, $R_{+,0,-}$ etc - readers are encouraged to have at it!} but a couple of points are worth noting in passing:\\
\\
\textbf{(i)} The range of the operator $R_{+,0,-}$ is a subspace of the space of periodic functions on $\mathbb{C}$ with period $1$ since it is immediate that for any remainder-summable function $f$ we have $R_{+,0,-}\left[f\right]\left(z+1\right)=R_{+,0,-}\left[f\right]\left(z\right)$. This allows us to make use of the theory of Fourier series. We will see an application of this in the next subsection when we use it to deduce the functional equation of the Gamma function within a generalised C\'{e}saro framework (and more generally see how it may be likewise applied to derive functional equations for other functions also).\\
\\
\textbf{(ii)} Calculations readily reveal that, for any $n\in\mathbb{Z}_{\geq{0}}$ we have $R_{+,0,-}\left[\tilde{z}^{n}\right]\left(z\right)=0$ for all $z\in\mathbb{C}$, i.e.
\begin{equation}
R_{+,0,-}\left[\tilde{z}^{n}\right]\equiv 0.
\label{bidirectional_z_to_n}\end{equation}
Thus non-negative integer powers of $z$, and hence all polynomials, lie in the kernel of the operator $R_{+,0,-}$ and it follows rapidly from this that the kernel of $R_{+,0,-}$ also contains all functions which consist of a product of a polynomial and a period-1 function. The fact that generically $R_{+,0,-}\left[\tilde{z}^{-s}\right]\neq 0$ for $s\notin\mathbb{Z}_{\leq{0}}$ suggests that it is the presence of a singularity in the input function which somehow leads to non-trivial output.\\
\\
\textbf{(iii)} Recalling that $R_{+}\left[\tilde{z}^{-s}\right](0)=\zeta(s)$ and that $(-n)^{-s}=e^{-i\pi s}n^{-s}$, equation \ref{bidirectional_z_to_n} allows us moreover to deduce that $\zeta(-2)=0=\zeta(-4)=\zeta(-6) \dots$, i.e. that zeta has trivial zeros at the negative even integer points. We omit details here, but this is a first example showing how considering bi-directional sums can often lead to interesting results, for the zeta function and for other problems.\\
\\
\textbf{Remainder products:} The concept of remainder summation allows us also to define a notion of remainder product via exponentiation of the remainder sum of the logarithm:
\begin{equation}
\prod_{R,+}[f](z):=\exp\left(R_{+}[\ln(f)](z)\right)
\label{Rem_prod}\end{equation}
Effectively, we are defining the remainder product of $f$, $\prod_{R,+}[f]$, by defining its logarithm as the remainder sum of $\ln(f)$. Note that this definition is agnostic to the choice of divergence-scheme used to make sense of the remainder sum on the RHS, but we shall only really be concerned with the case of using generalised geometric C\'{e}saro convergence in this and future papers. A variety of other possible approaches could also be used to define $\prod_{R,+}[f]$ - for example, we could take the partial-product function, $\prod_{j=1}^{k}f(z+j)$, and try to apply C\'{e}saro convergence directly to it - but the approach chosen is the best-adapted for most situations and the only one we will consider. Extensions of equation \ref{Rem_prod} to define $\prod_{R,+,0}[f]$, $\prod_{R,-}[f]$, $\prod_{R,+,0,-}[f]$ etc follow in the obvious way.\\
\\
\textbf{Extending finite products:} Just as we naturally extended the concept of a finite sum, $s(k)=\sum_{j=1}^{k}f(j)$, from a function of a discrete variable $k\in\mathbb{Z}_{>0}$ to a function, $s(z_{0})=\sum_{j=1}^{z_{0}}f(j)$, of arbitrary $z_{0}\in\mathbb{C}$ via equation \ref{Finite_sum_z_0}, we can also extend the traditional notion of a finite product. We extend $\pi(k)=\prod_{j=1}^{k}f(j)$, to a function of $z_{0}\in\mathbb{C}$ by writing $\prod_{j=1}^{z_{0}}f(j)=exp\left(\sum_{j=1}^{z_{0}}\ln\left(f(j)\right)\right)$ and using equation \ref{Finite_sum_z_0} to make sense of $\sum_{j=1}^{z_{0}}\ln\left(f(j)\right)$, i.e.
\begin{equation}
\prod_{j=1}^{z_{0}}f(j)\,=\,\frac{\prod_{R,+}[f](0)}{\prod_{R,+}[f](z_{0})}:=\exp\left(R_{+}[\ln(f)](0)-R_{+}[\ln(f)](z_{0})\right)
\label{Finite_prod}\end{equation}
\\
\textbf{Examples: (i)} If we take $f$ as the constant function $f(z)=e$ and recall from section 2.1, equation \ref{Hurewicz_zeta_0}, that $R_{+}\left[1\right](z_{0})=-z_{0}-\frac{1}{2}$, then we see at once that $\prod_{R,+}[e](z_{0})=e^{-z_{0}-\frac{1}{2}}$; and if we let $g(z_{0})$ denote the finite product $\prod_{j=1}^{z_{0}}e$, then $g(z_{0})=e^{z_{0}}$ (which fits with the fact that clearly $g(k)=e^{k}$ for $k\in\mathbb{Z}_{>0}$).\\
\\
\textbf{(ii)} Next consider the case of the identity function $f(z)=z$ and the finite product $\prod_{j=1}^{z_{0}}j$. Clearly, for $z_{0}=k\in\mathbb{Z}_{>0}$ this gives $k!$ so that, morally, the extension to all of $z_{0}\in\mathbb{C}$ should yield $\Gamma(z_{0}+1)$. As such we shall for now denote $\prod_{j=1}^{z_{0}}j$ by $\tilde{\Gamma}(z_{0}+1)$ and we note that, like $\Gamma$, it also clearly satisfies the recurrence relation that $\tilde{\Gamma}(z_{0}+1)=z_{0}\cdot\tilde{\Gamma}(z_{0})$ and has simple poles at $0$ and all negative integer points. It turns out that $\tilde{\Gamma}$ is in fact the same as $\Gamma$ - so that the remainder C\'{e}saro definition
\begin{equation}
\tilde{\Gamma}(z_{0}+1)\,=\,\frac{\prod_{R,+}[\tilde{z}](0)}{\prod_{R,+}[\tilde{z}](z_{0})}:=\exp\left(R_{+}[\ln(\tilde{z})](0)-R_{+}[\ln(\tilde{z})](z_{0})\right)
\label{Gamma_Cesaro_def}\end{equation}
actually gives the most natural possible way of extending the traditional factorial function to $\Gamma$ on all of $\mathbb{C}\setminus\mathbb{Z}_{<0}$. We devote the next section to three C\'{e}saro calculations which demonstrate this, and also derive its functional equation.

\subsection{The C\'{e}saro analysis of the Gamma function}

The three computations we perform are, first, to show that $\tilde{\Gamma}(z_{0}+1)$ is well-defined (i.e. that the RHS in equation \ref{Gamma_Cesaro_def} is remainder C\'{e}saro-summable to give a well-defined value for all $z_{0}\in\mathbb{C}\setminus\mathbb{Z}_{<0}$); secondly to derive the Taylor series for $\tilde{\Gamma}(z_{0}+1)$ around $z_{0}=0$ and hence deduce that $\tilde{\Gamma}$ is identical to $\Gamma$; and finally to show how the functional equation for $\Gamma$ drops out of this C\'{e}saro framework in a natural way by combining bi-directional summation and Fourier theory, as foreshadowed in the previous section.\\
\\
\textbf{Calculation 1 - $\tilde{\Gamma}(z_{0}+1)$ is well-defined for all $z_{0}\in\mathbb{C}\setminus\mathbb{Z}_{<0}$:} For any given $z_{0}\notin\mathbb{Z}_{<0}$, let $z=z_{0}+k+\alpha$ in the usual way. Applying the Euler-McLaurin sum formula to calculate the p-sum function for $R_{+}\left[\ln \tilde{z}\right](z_{0})$, we have
\begin{eqnarray*}
\sum_{j=1}^{k}\ln(z_{0}+j) & = & \intop^{k}\ln(z_{0}+t)\,\textrm{d}t+\frac{1}{2}\ln(z_{0}+k)+C_{z_{0}}+o(1)\\
 & = & (z_{0}+k+\frac{1}{2})\ln(z_{0}+k)-(z_{0}+k)+C_{z_{0}}+o(1)\\
 & \overset{C}{\sim} & (z_{0}+k+\alpha)\ln(z_{0}+k+\alpha)-(z_{0}+k+\alpha)+C_{z_{0}}\\
 & \overset{C}{\sim} & C_{z_{0}} \qquad .
\end{eqnarray*}
Here we are noting that $\ln(z_{0}+k+\alpha)=\ln(z_{0}+k)+\frac{\alpha}{z_{0}+k}+\dots$ and that $(\alpha-\frac{1}{2})\ln(z_{0}+k)$ converges to $0$ under a single application of the C\'{e}saro averaging operator $P$; and also recalling that $z\ln z$ and $z$ are generalised eigenfunctions of $P$ in the geometric variable $z$, with eigenvalue $\frac{1}{2}$. Thus $R_{+}\left[\ln \tilde{z}\right](z_{0})$ is well-defined (as $C_{z_{0}}$) for arbitrary such $z_{0}$ and so in equation \ref{Gamma_Cesaro_def} the RHS makes sense, giving $\tilde{\Gamma}$ as the well-defined function
\begin{equation}
\tilde{\Gamma}(z_{0}+1)=\exp(C_{0}-C_{z_{0}})=\frac{\exp(C_{0})}{\exp(C_{z_{0}})}
\label{Remprodza}\end{equation}
on all of $\mathbb{C\setminus}\{\mathbb{Z}_{<0}\}$, where
\begin{equation}
C_{z_{0}}=\underset{k\rightarrow\infty}{\lim}\left\{ \sum_{j=1}^{k}\ln(z_{0}+j)-(z_{0}+k+\frac{1}{2})\ln(z_{0}+k)+(z_{0}+k)\right\}
\label{Remprodzb}\end{equation}\\
\\
\textbf{Comments on calculation 1:} \textbf{(i)} It is \textit{critical} in the above generalised C\'{e}saro analysis that summands are introduced geometrically at $z_{0}+j$ and that the C\'{e}saro eigenfunctions we then annihilate in taking C\'{e}saro limits are functions of the geometric variable $z=z_{0}+k+\alpha$, rather than simply the summation index $k$ or the contour parameter $t=k+\alpha$. Indeed, it is easy to see that if we used $t$ rather than $z$ as our C\'{e}saro parameter, then we would always end up with residual pure log-divergences in our p-sum function above - so that we would not be able to assign a C\'{e}saro limit, and hence a value for $\tilde{\Gamma}(z_{0}+1)$, at even a single point $z_{0}\in\mathbb{C}$. This is why we have emphasised in earlier sections that generalised C\'{e}saro summation/convergence is fundamentally \textit{geometric} in nature, and sensitive to the locations of summands (and in future papers, of roots, poles etc).\\
\\
\textbf{(ii)} In this regard note also that, unlike the case of $\zeta$ where the defining sum was already classically convergent on the region $Re(s)>1$ and generalised C\'{e}saro methodology was only needed to extend to the region of divergence $Re(s)\leq 1$, this is a situation where C\'{e}saro methods are required to make sense of $\tilde{\Gamma}(z_{0}+1)$ at all $z_{0}\in\mathbb{C}\setminus\mathbb{Z}_{<0}$\footnote{Other than the points $z_{0}=k\in\mathbb{Z}_{\geq{0}}$, of course, where the defining product trivially collapses to $k!$ }.\\
\\
\textbf{(iii)} When $z_{0}=0$ the value of $C_{0}$ in equations \ref{Remprodza} and \ref{Remprodzb} corresponds to $-\zeta^{\prime}(0)$ (which has the same defining \textit{picture} as $\tilde{\Gamma}^{\prime}(1)$ since formally $\frac{d}{ds}\left(\sum_{n=1}^{\infty}n^{-s}\right)\big|_{s=0}=-\sum_{n=1}^{\infty}\ln(n)$ and the summands are likewise introduced at the points $n=1,2,3,\dots$). It is known that $\zeta^{\prime}(0)=-\frac{1}{2}\ln(2\pi)$ (although we omit the digression required to give a C\'{e}saro derivation of this fact here) and so the definition of $\tilde{\Gamma}(z_{0}+1)$ in equation \ref{Remprodza} can be rewritten more simply as just
\begin{equation}
\tilde{\Gamma}(z_{0}+1)\,=\,\frac{\left(2\pi\right)^{\frac{1}{2}}}{\prod_{R,+}\left[\tilde{z}\right](z_{0})}\,=\,\exp\left(\frac{1}{2}\ln(2\pi)-R_{+}\left[\ln(\tilde{z})\right](z_{0})\right)
\label{Gamma_tilde_formula}\end{equation}
\\
\textbf{Calculation 2 - Calculation of Taylor series of $\tilde{\Gamma}(z_{0}+1)$ near $0$ and proof that $\tilde{\Gamma}=\Gamma$:} We now turn to the Taylor series near $0$ of $\tilde{\Gamma}(z_{0}+1)$, or rather of $\ln(\tilde{\Gamma}(z_{0}+1))$, since equation \ref{Gamma_tilde_formula} allows us to write
\begin{equation}
\ln(\tilde{\Gamma}(z_{0}+1))\,=\,\frac{1}{2}\ln(2\pi)-R_{+}\left[\ln(\tilde{z})\right](z_{0})
\label{ln_Gamma_tilde_formula}\end{equation}
and we can thus perform C\'{e}saro calculations directly for $\ln(\tilde{\Gamma}(z_{0}+1))$ without also having to worry about exponentiation\footnote{Peeking ahead to the fact that $\tilde{\Gamma}=\Gamma$, this itself explains in C\'{e}saro terms the long-standing observation of mathematicians that $\ln\left(\Gamma(z_{0}+1)\right)$ is often easier to analyse than $\Gamma(z_{0}+1)$ itself}.

Clearly at $z_{0}=0$ we have $\ln(\tilde{\Gamma}(1))=0$ and for $\epsilon$ small we have that formally
\begin{equation*}
\ln(\tilde{\Gamma}(\epsilon+1))\,=\,\frac{1}{2}\ln(2\pi)-R_{+}\left[\ln(\tilde{z})\right](\epsilon)\,=\,\frac{1}{2}\ln(2\pi)-\sum_{j=1}^{\infty}\ln(j+\epsilon)
\end{equation*}
while termwise for each $j$ we have
\begin{equation*}
\ln(j+\epsilon)\,=\,\ln j + \ln\left(1+\frac{\epsilon}{j}\right)\,=\,\ln j + \frac{\epsilon}{j}-\frac{1}{2}\frac{\epsilon^{2}}{j^{2}}+\frac{1}{3}\frac{\epsilon^{3}}{j^{3}}-\dots \qquad .
\end{equation*}
For $n\in\mathbb{Z}_{\geq{2}}$ the coefficient of $\epsilon^{n}$ is the convergent sum $\frac{(-1)^{n}}{n}\sum_{j=1}^{\infty}\frac{1}{j^{n}}\,=\,\frac{(-1)^{n}}{n}\zeta(n)$. 

However, the same naive calculation for the coefficient of $\epsilon^{1}$ in the Taylor series gives the negative of the divergent sum $\sum_{j=1}^{\infty}\frac{1}{j}$. It would therefore seem at first glance that this divergence cannot be handled even within a C\'{e}saro framework because the p-sum function appears to be $s(z)=s(k+\alpha)=-\sum_{j=1}^{k}\frac{1}{j}=-\ln k - \gamma+o(1)=-\ln z - \gamma+o(1)$ and the pure log-divergence, $\ln z$, cannot be assigned a generalised C\'{e}saro limit. 

Fortunately, however, a more careful C\'{e}saro calculation - specifically one that again pays critical attention to the "picture" involved and adds in summands with due care for their geometric placement - resolves this apparent problem and leaves us with the Taylor series coefficient of $\epsilon^{1}$ as minus the Euler-Mascheroni constant $\gamma\approx 0.577$.

To see this, recall that the coefficient of $\epsilon^{1}$ is given by 
\begin{eqnarray*}
\frac{d}{dz_{0}}\left(\ln\left(\tilde{\Gamma}(z_{0}+1)\right)\right)\Big|_{z_{0}=0} & = & \lim_{h\rightarrow0}\;\frac{\ln(\tilde{\Gamma}(1+h))-\ln(\tilde{\Gamma}(1))}{h}\\
 & = & \lim_{h\rightarrow0}\;\frac{R_{+}\left[\ln(\tilde{z})\right](0)-R_{+}\left[\ln(\tilde{z})\right](h)}{h}
\end{eqnarray*}
The term $R_{+}\left[\ln(\tilde{z})\right](0)$ gives rise to a set of terms, $\ln j$, each added at the point $j$ on the real axis. The term $R_{+}\left[\ln(\tilde{z})\right](h)$ gives rise to a set of terms, $\ln(j+h)$, and the error in the above calculation arises from forgetting the proper location of these summands and treating them as though they too are subtracted at point $j$, leaving the residual sum of $\sum_{j=1}^{\infty}\frac{1}{j}$ (after taking the limit as $h\rightarrow 0$), with its untreatable log-divergence.

In fact, however, each term $\ln(j+h)$ should be subtracted not at $j$ but at $j+h$ and this, as Prof Robert Frost remarked (pers. comm.), makes all the difference. 

When we do this, and ignore overall $O(h)$ components, the p-sum function for calculation of $\frac{d}{dz_{0}}\left(\ln\left(\tilde{\Gamma}(z_{0}+1)\right)\right)\Big|_{z_{0}=0}$ becomes as pictured in Fig. 1.

\includegraphics[height=8.5cm, width=10cm]{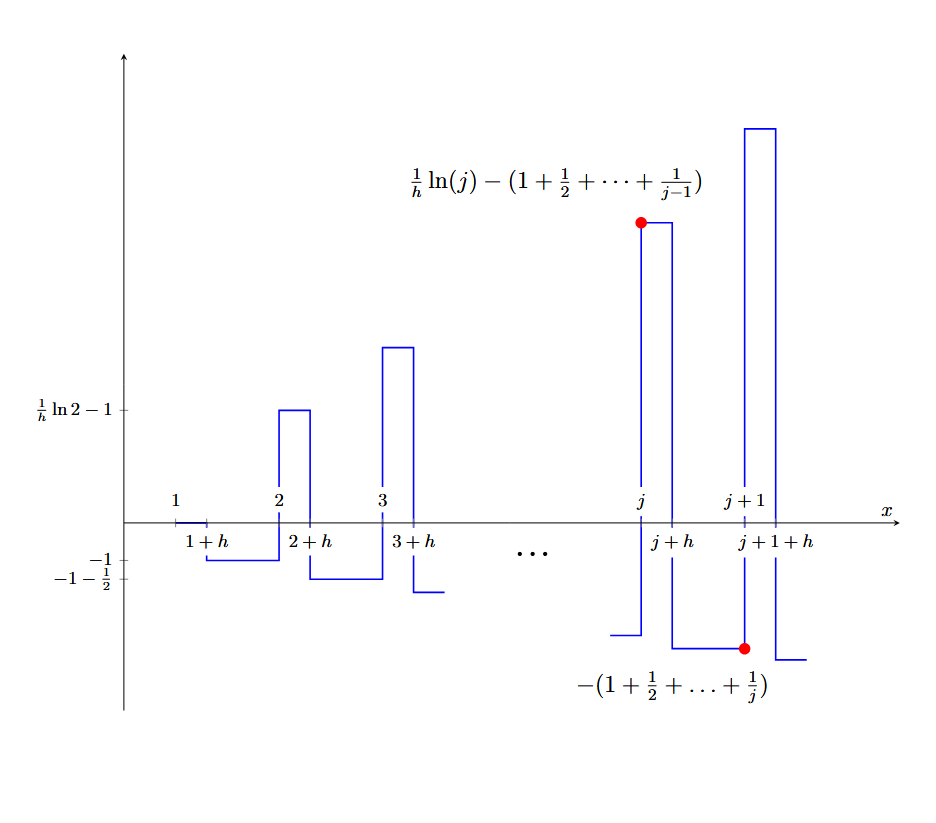}
\\
This p-sum function is now a quasi-oscillatory divergent function on $[0,\infty)$ whose value on a generic interval $[j,j+1)$ is given by
\begin{equation*}
s(z)=\begin{cases}
\frac{1}{h}\ln(j) - \{1+\frac{1}{2}+\frac{1}{3}+\cdots+\frac{1}{j-1}\}\quad, & j\leq{z}<j+h\\
-\{1+\frac{1}{2}+\frac{1}{3}+\cdots+\frac{1}{j}\}\quad\quad, & j+h<z<j+1\end{cases}
\end{equation*}
Since the integral of this p-sum function over $[j,j+1)$ gives
\begin{equation*}
\int_{j}^{j+1}s(z)\,dz\,=\,\ln(j)-\sum_{l=1}^{j}\frac{1}{l}+O(h)\,=\,-\gamma+O\left(\frac{1}{j}\right)+O(h)
\end{equation*}
it follows that the average value of $s(z)$ in C\'{e}saro terms is $-\gamma$, i.e. $P[s](z)\rightarrow -\gamma$ classically as $z\rightarrow\infty$. Thus we have 
\begin{equation}
\frac{d}{dz_{0}}\left(\ln\left(\tilde{\Gamma}(z_{0}+1)\right)\right)\Big|_{z_{0}=0}=-\gamma
\label{Gamma_Tilde_First_Deriv}\end{equation} 
as claimed.\\
\\
\textbf{Comment:} Not only does this calculation of $\frac{d}{dz_{0}}\left(\ln\left(\tilde{\Gamma}(z_{0}+1)\right)\right)\Big|_{z_{0}=0}$ once again demonstrate the criticality of respecting geometry within the remainder C\'{e}saro framework, it accomplishes what in a Physics context would be understood as a renormalisation calculation of the divergent series $-\sum_{j=1}^{\infty}\frac{1}{j}$ as $-\gamma$.\\
\\
Combining the above calculations of the coefficients of $\epsilon^{n}$ for $n\in\mathbb{Z}_{\geq{2}}$ and for $n=1$, we have that the Taylor series for $\ln\left(\tilde{\Gamma}(z_{0}+1)\right)$ around $z_{0}=0$ is
\begin{equation}
\ln\left(\tilde{\Gamma}(z_{0}+1)\right)\,=\,-\gamma z_{0}+\sum_{n=2}^{\infty}\frac{(-1)^{n}}{n}\,\zeta(n)\,z_{0}^{n}
\label{Gamma_tilde_TS}\end{equation}
with radius of convergence $1$. Since this agrees with the known Taylor series for $\Gamma(z_{0}+1)$ on the unit disc, it follows that $\tilde{\Gamma}(z_{0}+1)=\Gamma(z_{0}+1)$ on this disc and hence, by the uniqueness of analytic continuation, on all of $\mathbb{C}\setminus\mathbb{Z}_{<0}$. 

Thus our earlier moral conviction - that $\tilde{\Gamma}$ defined as $\tilde{\Gamma}(z+1)=\prod_{j=1}^{z}j$ and understood by equation \ref{Gamma_Cesaro_def} should simply be the classical $\Gamma$-function - is vindicated. On this basis we omit the distinction between them and simply write $\tilde{\Gamma}$ as $\Gamma$ going forward.\\
\\
\textbf{Calculation 3 - C\'{e}saro derivation of the functional equation of $\Gamma$:} Having shown that our generalised C\'{e}saro definition in equation \ref{Gamma_Cesaro_def} does indeed yield $\Gamma(z_{0}+1)$, we show that the famous functional equation for $\Gamma$, namely that
\begin{equation}
\Gamma(z_{0})\,\Gamma(1-z_{0})\,=\,\frac{\pi}{\sin\left(\pi z_{0}\right)}
\label{GammaFunctionalEquation_2}\end{equation}
can be deduced in a systematic way within this C\'{e}saro framework.

Taking logarithms, equation \ref{GammaFunctionalEquation_2} is equivalent to saying that
\begin{equation}
\ln(\Gamma(z_{0}))+\ln(\Gamma(1-z_{0}))\,=\,\ln(\pi) - \ln(\sin\left(\pi z_{0}\right))
\label{LogGammaFunctionalEquation}\end{equation}
and by equation \ref{ln_Gamma_tilde_formula} this is equivalent to proving that
\begin{equation}
R_{+,0}\left[\ln\tilde{z}\right](z_{0})+R_{+,0}\left[\ln\tilde{z}\right](1-z_{0})\,=\,\ln(2)+\ln(\sin\left(\pi z_{0}\right))
\label{LogGammaFunctionalEquation_2}\end{equation}
where we have noted that the LHS in equation \ref{ln_Gamma_tilde_formula} can be changed from $\ln(\Gamma(z_{0}+1))$ to $\ln(\Gamma(z_{0}))$ as long as we change the operator on the RHS from $R_{+}$ to $R_{+,0}$.

\includegraphics[height=9cm, width=10cm]{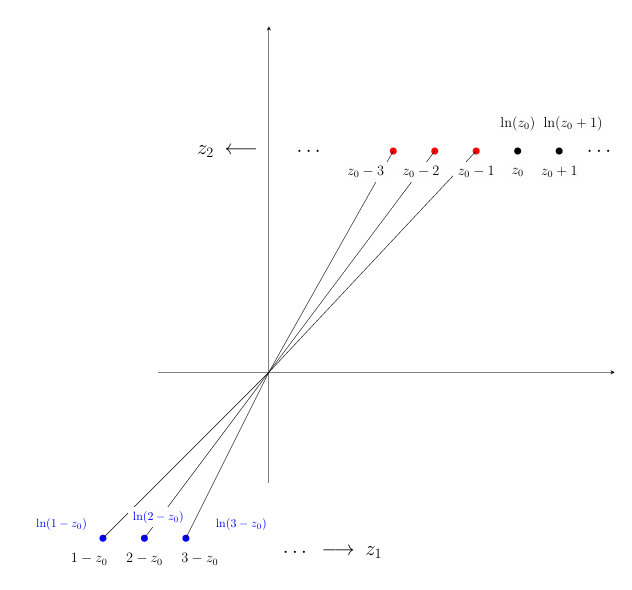}
\\
Now consider figure 2. The term $R_{+,0}\left[\ln\tilde{z}\right](z_{0})$ in equation \ref{LogGammaFunctionalEquation_2} introduces $\ln\tilde{z}$-summands at the points $z_{0}, z_{0}+1, z_{0}+2, \ldots $ shown in black and the term $R_{+,0}\left[\ln\tilde{z}\right](1-z_{0})$ introduces $\ln\tilde{z}$-summands at the points $-z_{0}+1, -z_{0}+2, -z_{0}+3, \ldots $ shown in blue. As shown in red in the figure, these latter points are the negatives of the set $z_{0}-1, z_{0}-2, z_{0}-3, \ldots $. 

Since we can easily relate the values of $\ln w$ and $\ln(-w)$, we can thus try to replace $R_{+,0}\left[\ln\tilde{z}\right](1-z_{0})$ in equation \ref{LogGammaFunctionalEquation_2} with $R_{-}\left[\ln\tilde{z}\right](z_{0})$, albeit noting that we will need to show appropriate C\'{e}saro care in doing so in order to reflect the fact that the \textit{geometric} variable, $z$, whose eigenfunctions/generalised eigenfunctions we annihilate is different in each case - for $R_{+,0}\left[\ln\tilde{z}\right](1-z_{0})$ we will have $z=z_{1}:=-z_{0}+k+\alpha$, while for $R_{-}\left[\ln\tilde{z}\right](z_{0})$ we will have $z=z_{2}:=z_{0}-k-\alpha=-z_{1}$. 

Taken together with the $R_{+,0}\left[\ln\tilde{z}\right](z_{0})$ term, this $R_{-}\left[\ln\tilde{z}\right](z_{0})$ term then gives us overall $R_{+,0,-}\left[\ln\tilde{z}\right](z_{0})$ on the LHS in equation \ref{LogGammaFunctionalEquation_2} and since, as discussed earlier in section 2.2, $R_{+,0,-}\left[\ln\tilde{z}\right](z_{0})$ is automatically a period-1 function, we can then apply Fourier theory to it to deduce its equivalence to the RHS.

Let us carry out this cunning plan in detail. We take $-\pi<\theta\leq{\pi}$ as the principal branch of $\ln$ and consider first the case of $Im(z_{0})>0$ (the case of $Im(z_{0})<0$ follows in the same way and we leave this and the further minor extension to the real axis to the reader).

In this case we have $\ln(-z_{0}+j)=\ln(z_{0}-j)-i\pi$ for all $j\in\mathbb{Z}$. Using the geometric variables $z_{1}$ and $z_{2}$ just defined, we denote the p-sums for $R_{+,0}\left[\ln(\tilde{z})\right](1-z_{0})$ and $R_{-}\left[\ln(\tilde{z})\right](z_{0})$ by $s_{1}$ and $s_{2}$ respectively, so that
\begin{equation}
R_{+,0}[\ln\tilde{z}](1-z_{0})=\underset{z_{1}\rightarrow\infty}{Clim}\, s_{1}(z_{1})\quad\textrm{and}\quad R_{-}[\ln\tilde{z}](z_{0})=\underset{z_{2}\rightarrow-\infty}{Clim}\, s_{2}(z_{2}) \quad .
\label{R+0_s1_R-_s2_relationship}\end{equation}
Then we have that
\begin{eqnarray*}
s_{1}(z_{1}) & = & \sum_{j=1}^{k}\ln(j-z_{0})=\sum_{j=1}^{k}\ln(z_{0}-j)-i\pi\sum_{j=1}^{k}1\\
 & = & s_{2}(z_{2})-i\pi k
\end{eqnarray*}
Now as before in section 2.3, we have that
\begin{equation*}
s_{1}(z_{1})=(-z_{0}+k+\frac{1}{2})\ln(-z_{0}+k)-(-z_{0}+k)+C_{-z_{0}}+o(1)\overset{C}{\sim}C_{-z_{0}}
\end{equation*}
Therefore
\begin{equation*}
s_{2}(z_{2})=-(z_{0}-k-\frac{1}{2})[\ln(z_{0}-k)-i\pi]+(z_{0}-k)+C_{-z_{0}}+i\pi k+o(1).
\end{equation*}
But, via the same sort of elementary C\'{e}saro calculations we have already performed several times in this paper, we have that
\begin{eqnarray*}
-z_{2}\ln z_{2} & = & -(z_{0}-k-\alpha)\ln(z_{0}-k-\alpha)\\
 & = & -(z_{0}-k)\cdot\{\ln(z_{0}-k)-\frac{\alpha}{(z_{0}-k)}\}+\alpha\ln(z_{0}-k)+o(1)\\
 & \overset{c}{\sim} & -(z_{0}-k)\ln(z_{0}-k)+\frac{1}{2}\ln(z_{0}-k)+\frac{1}{2} \quad .
\end{eqnarray*}
Thus
\begin{eqnarray*}
s_{2}(z_{2}) & \overset{C}{\sim} & -z_{2}\ln z_{2}+z_{2}+C_{-z_{0}}+i\pi(z_{0}-\frac{1}{2})\\
 & \overset{C}{\sim} & C_{-z_{0}}+i\pi(z_{0}-\frac{1}{2})
\end{eqnarray*}
and it follows at once from equation \ref{R+0_s1_R-_s2_relationship} that, as hoped, we do have a relationship between $R_{+,0}\left[\ln(\tilde{z})\right](1-z_{0})$ and $R_{-}\left[\ln(\tilde{z})\right](z_{0})$, namely
\begin{equation}
R_{+,0}\left[\ln(\tilde{z})\right](1-z_{0})\,=\,R_{-}\left[\ln(\tilde{z})\right](z_{0})-i\pi(z_{0}-\frac{1}{2}) \quad .
\label{R+0_R-_relationship}\end{equation}
In equation \ref{LogGammaFunctionalEquation_2} it follows that proving the functional equation for $\Gamma$ is equivalent to proving that
\begin{equation}
R_{+,0,-}\left[\ln(\tilde{z})\right](z_{0})\,=\,i\pi(z_{0}-\frac{1}{2})+\ln(2)+\ln(\sin\left(\pi z_{0}\right)) \quad .
\label{LogGammaFunctionalEquation_3}\end{equation}
Now since, as noted, $R_{+,0,-}\left[\ln(\tilde{z})\right](z_{0})$ is a periodic function of $z_{0}$ with period $1$ it follows that it has the Fourier series expansion
\begin{equation}
R_{+,0,-}\left[\ln(\tilde{z})\right](z_{0})\,=\,\sum_{n=-\infty}^{\infty}a_{n}e^{2\pi inz_{0}}
\label{R+0-_FS_1}\end{equation}
where
\begin{eqnarray*}
a_{n} & = & \int_{0}^{1}R_{+,0,-}\left[\ln\tilde{z}\right](z_{0})e^{-2\pi inz_{0}}\,dz_{0}\,=\,\int_{-\infty}^{\infty}\ln(z_{0})e^{-2\pi inz_{0}}\,dz_{0}\\
 & = & \mathcal{F}\left[\ln\tilde{z}\right](n) \quad .
\end{eqnarray*}
Here $\mathcal{F}$ is the Fourier transform defined by
\begin{equation}
\mathcal{F}\left[f\right](\xi)\,:=\,\int_{-\infty}^{\infty}f(x)e^{-2\pi ix\xi}\,dx \quad .
\label{FT_definition}\end{equation}
Now  on $\mathbb{R}$
\begin{equation*}
\ln \tilde{z} = \ln |\tilde{z}| + i\pi (1-u(\tilde{z}))
\end{equation*}
where $u(\tilde{z})$ is the one-sided Heaviside function given by
\begin{equation*}
u(\tilde{z})=\begin{cases}
0 \quad, & \tilde{z}<0\\
1 \quad, & \tilde{z}>0\end{cases}
\end{equation*}
and the Fourier transforms of all these component-pieces can be looked up from standard tables as
\begin{equation*}
\mathcal{F}\left[\ln |\tilde{z}|\right](\xi) = -\frac{1}{2}\frac{1}{|\xi|}+\gamma\delta_{0}(\xi)\; ,
\end{equation*}
\begin{equation*}
\mathcal{F}\left[1\right](\xi) = \delta_{0}(\xi)
\end{equation*}
and
\begin{equation*}
\mathcal{F}\left[u(\tilde{z})\right](\xi) = \frac{1}{2\pi i}\frac{1}{\xi}+\frac{1}{2}\delta_{0}(\xi)\; .
\end{equation*}
It follows at once that for $n\in\mathbb{Z}_{>0}$ we have $a_{n}=-\frac{1}{n}$, while for $n\in\mathbb{Z}_{<0}$ we have $a_{n}=0$ and finally $a_{0}=0$ also, by an elementary C\'{e}saro calculation. In equation \ref{LogGammaFunctionalEquation_3} we therefore have
\begin{equation*}
LHS = -\sum_{n=1}^{\infty}\frac{e^{2\pi inz_{0}}}{n} \; .
\end{equation*}
But, on the other hand, since $Im(z_{0})>0$ it follows that $|e^{2\pi inz_{0}}|<1$ and so
\begin{equation*}
\ln(\sin(\pi z_{0}))=\ln\left(\frac{-e^{-i\pi z_{0}}}{2i}\left(1-e^{2\pi iz_{0}}\right)\right)=-i\pi z_{0}+\frac{i\pi}{2}-\ln2-\sum_{n=1}^{\infty}\frac{e^{2\pi inz_{0}}}{n}
\end{equation*}
from which it follows that the RHS in equation \ref{LogGammaFunctionalEquation_3} is also given by
\begin{equation*}
RHS = -\sum_{n=1}^{\infty}\frac{e^{2\pi inz_{0}}}{n} \; .
\end{equation*}
Hence equation \ref{LogGammaFunctionalEquation_3} is verified and this completes the proof of the functional equation of $\Gamma$, equation \ref{GammaFunctionalEquation_2}.\\
\\
\textbf{Comments:} \textbf{(i)} The plan outlined and followed in the above proof should potentially work not just for $\ln(\Gamma)$ (and hence $\Gamma$ after exponentiation) but for any function, $g(z_{0}):=R_{+,0}[f](z_{0})$ defined by a remainder C\'{e}saro sum. As long as we can relate $f(-w)$ to $f(w)$ then, in light of figure 2, we can potentially convert $R_{+,0}[f](z_{0})+R_{+,0}[f](1-z_{0})$ into an expression involving $R_{+,0,-}[f](z_{0})$ and then utilise Fourier theory to understand this bi-directional sum and obtain a functional equation that relates $g(z_{0})$ and $g(1-z_{0})$. C\'{e}saro remainder summation thus gives one general avenue by which functional equations can sometimes be derived - and one which shows why many such functional equations exhibit a common structure in relating values of the function at $z_{0}$ and $1-z_{0}$.\\
\\
\textbf{(ii)} For example, in the case of $\zeta_{H}(z_{0};s)=R_{+}\left[\tilde{z}^{-s}\right](z_{0})$ where we have $f(\tilde{z})=\tilde{z}^{-s}$ we can enact the same program and thereby end up with a functional equation for $\zeta_{H}(z_{0};s)$ that relates its values at $z_{0}$ and $1-z_{0}$ for any s. It is not as elegant a relationship as for $\ln(\Gamma(z_{0}+1))$ because the Fourier series that arises does not simplify in general. But if we specialise to $z_{0}=\frac{1}{2}$ then we \textit{can} simplify materially and in fact we then obtain the functional equation for the Riemann zeta function, $\zeta(s)$, as the residual relationship in $s$ only. We shall, however, defer the details of this derivation to Appendix 4.1.

\section{Dilation and scaling invariance of C\'{e}saro convergence}

We turn now to the dilation and scaling invariance properties of generalised geometric C\'{e}saro convergence, starting with dilation-invariance.

The Dilation group is a group of operators $\lbrace D_{r}\rbrace_{r\in\mathbb{R}_{>0}}$ under multiplication (so $D_{r}\cdot D_{s}=D_{rs}$) given by $D_{r}[f](z):=f(rz)$. Our key result is that\\
\\
\textbf{Theorem 1:} \textit{Generalised geometric C\'{e}saro convergence is dilation-invariant, i.e. for any} $r\in\mathbb{R}_{>0}$\textit{, if} $g(z):=D_{r}[f](z)=f(rz)$ \textit{then} $g$ \textit{has a generalised C\'{e}saro limit along contour $\gamma$ if and only if} $f$ \textit{has a generalised C\'{e}saro limit along contour $\tilde{\gamma}:=r\gamma$ and in this case} $\underset{z\rightarrow\infty}{Clim}\,g(z)\,=\,\underset{z\rightarrow\infty}{Clim}\,f(z)$.\\
\\
We give two proofs, the first direct, the second more abstract.\\
\\
\textbf{Proof 1:} By definition we have
\begin{eqnarray*}
P\left[g\right](\gamma(X)) & = & \frac{1}{X}\int_{0}^{X}f(rx)\,dx\\
 & = & \frac{1}{rX}\int_{0}^{rX}f(u)\,du \,=\, P\left[f\right](\tilde{\gamma}(rX))
\end{eqnarray*}
Thus if $f\overset{C}{\rightarrow}L$ along $\tilde{\gamma}$ under $P$ then $g\overset{C}{\rightarrow}L$ along $\gamma$ under $P$ also, and vice-versa; and clearly this continues to hold for higher powers of $P$ in the same way. Furthermore, it is clear that:\\
\\
(i) Each one-dimensional eigenspace of $P$ is invariant under dilation (since $(rz)^{\rho}=r^{\rho}z^{\rho}$), and that\\
\\
(ii) The span of the space of generalised eigenfunctions of $P$ for any given eigenvalue is invariant under dilation (since $\ln(rz)=\ln(z)+\ln(r)$ so that for any $m\in\mathbb{Z}_{>0}$, $(rz)^{\rho}(\ln(rz))^{m}$ is a linear combination of generalised eigenfunctions $z^{\rho}(\ln z)^{n}$ with $0\leq{n}\leq{m}$).\\
\\
The theorem follows at once from definition 1.\\

\textbf{Proof 2:} The second proof uses the notion of the \textit{generator} of a group of operators like $\lbrace D_{r}\rbrace_{r\in\mathbb{R}_{>0}}$, or like the group of translations $\lbrace T_{h}\rbrace_{h\in\mathbb{R}}$ given by $T_{h}[f](z):=f(z+h)$. Recall that the generator of translations is $\frac{d}{dz}$ because for $h=\epsilon$ small $T_{\epsilon}[f](z_{0})=f(z_{0}+\epsilon)=f(z_{0})+\epsilon f^{\prime}(z_{0})+O(\epsilon^{2})$; and that as a consequence we have that $T_{h}=\exp(h\frac{d}{dz})$, which is simply a compact way of expressing Taylor's theorem:
\begin{equation*}
f(z_{0}+h)\,=\,\sum_{n=0}^{\infty}f^{(n)}(z_{0})\cdot \frac{h^{n}}{n!}\,=\,f(z_{0})+f^{\prime}(z_{0})\cdot h+\frac{1}{2!}f^{\prime\prime}(z_{0})\cdot h^{2}+\ldots \quad .
\end{equation*}
In the same way, after adjusting for the fact that $\lbrace D_{r}\rbrace$ is a group on $\mathbb{R}_{>0}$ under multiplication, rather than a group on $\mathbb{R}$ under addition (so that we need to consider $D_{\exp(\epsilon)}$ rather than $D_{\epsilon}$) we have that, working to order $\epsilon$,
\begin{equation*}
f(\exp(\epsilon)z_{0})\,=\,f((1+\epsilon)z_{0})\,=\,f(z_{0}+\epsilon z_{0})\,=\,f(z_{0})+\epsilon z_{0}\cdot f^{\prime}(z_{0})+O(\epsilon^{2}) \, .
\end{equation*}
It follows that the generator of dilations is
\begin{equation}
H_{D}=z\frac{d}{dz}\,.
\label{Generator_dilations}\end{equation}
On the other hand, from the definition of $P$ that $P[f](z_{0})=\frac{1}{z_{0}}\int_{0}^{z_{0}}f(z)\,dz$ it is elementary that $P^{-1}$ is given by
\begin{equation}
P^{-1}\,=\,\frac{d}{dz}\circ z\,=\,z\frac{d}{dz}+1\,=\,H_{D}+1 \quad .
\label{P_-1_def}\end{equation}
It follows that $P^{-1}$, and therefore obviously also $P$, commutes with the generator of dilations and hence with all dilations $D_{r}=\exp(\ln(r)H_{D})$; and hence that any regular polynomial $q(P)$ will likewise commute with all dilations, from which the theorem follows at once.\\
\\
\textbf{Comment:} Note that $z\frac{d}{dz}=\frac{d}{d\ln z}$ and that this observation and this operator turns up frequently in many aspects of the analysis of $\zeta$ (see e.g. [3] and its treatment in chapters 10 and 11 of Hardy's proof that there are infinitely many non-trivial Riemann zeros on the critical line, a topic to which we shall return in a number of future papers). The fact that $P^{-1}$ is almost identical to this operator probably explains why, as we have seen in this paper and its predecessor, the generalised \textit{continuous} C\'{e}saro framework is so well-adapted to the analysis of $\zeta$ and similar Dirichlet series, and why it works more smoothly than the corresponding \textit{discrete} C\'{e}saro framework for this purpose.

We shall consider also the invariance of generalised geometric C\'{e}saro convergence under another operator group, the scaling group. But first we devote the next subsection to two calculations (and one corollary result) demonstrating that C\'{e}saro dilation-invariance is not merely an abstract curiosity, but rather a powerful theorem that leads to interesting and non-trivial results.

\subsection{Three examples of C\'{e}saro dilation-invariance}

The calculations we undertake are the derivation of the multiplication formulae for $\Gamma(z_{0}+1)$; the deduction, as a corollary of this calculational method, of a corresponding result for $\zeta_{H}(z_{0};s)$; and the proof of the generalised integral result for $\zeta_{H}$ given earlier in equation \ref{zeta_Integral_Identity_3}.\\
\\
\textbf{Calculation 1: Derivation of the multiplication formulae for $\Gamma(z_{0}+1)$:} The famous multiplication formulae for $\Gamma$ state that for any $n\in\mathbb{Z}_{\geq{1}}$
\begin{equation}
(2\pi)^{\frac{(n-1)}{2}}\Gamma(z_{0}+1)=n^{z_{0}+\frac{1}{2}}\Gamma(\frac{z_{0}+1}{n})\cdot\Gamma(\frac{z_{0}+2}{n})\cdot\ldots\cdot\Gamma(\frac{z_{0}+n}{n}) \quad .
\label{GammaMultiplication}\end{equation}
On taking logs as usual and trivially re-arranging, this is equivalent to saying that
\begin{equation}
\ln(\Gamma(z_{0}+1))=\frac{1}{2}\ln(2\pi)+(z_{0}+\frac{1}{2})\ln n - \frac{n}{2}\ln(2\pi)+\sum_{l=1}^{n}\ln(\Gamma(\frac{z_{0}}{n}-\frac{(n-l)}{n}+1))
\label{GammaMultiplication_2}\end{equation}
and in light of equation \ref{ln_Gamma_tilde_formula}, this is in turn equivalent to the claim that
\begin{equation}
R_{+}[\ln \tilde{z}](z_{0})=(-z_{0}-\frac{1}{2})\ln n +\sum_{l=1}^{n}R_{+}[\ln \tilde{z}](\frac{z_{0}}{n}-\frac{(n-l)}{n}) \quad .
\label{GammaMultiplication_3}\end{equation}

But in light of figure 3 below, this result follows almost trivially from the dilation-invariance of C\'{e}saro convergence and our earlier calculation giving that $R_{+}[1](z_{0}) = -z_{0}-\frac{1}{2}$.

If we take the remainder C\'{e}saro sum for each $l$ on the RHS in turn and dilate it by $n$ (i.e. place each summand $\ln(\frac{z_{0}}{n}-\frac{(n-l)}{n}+j)$ at $z_{0}-(n-l)+jn$); and if we further note that $\ln(\frac{z_{0}-(n-l)+jn}{n})=\ln(z_{0}-(n-l)+jn)-\ln n$, then, taken over all the $l=1,2,\ldots,n$, the summands now intersperse perfectly to exactly cover the points $z_{0}+1,z_{0}+2,z_{0}+3,\ldots$ with each point $z_{0}+m$ having the summand $\ln(z_{0}+m)-\ln(n)$.

\includegraphics[height=9cm, width=10cm]{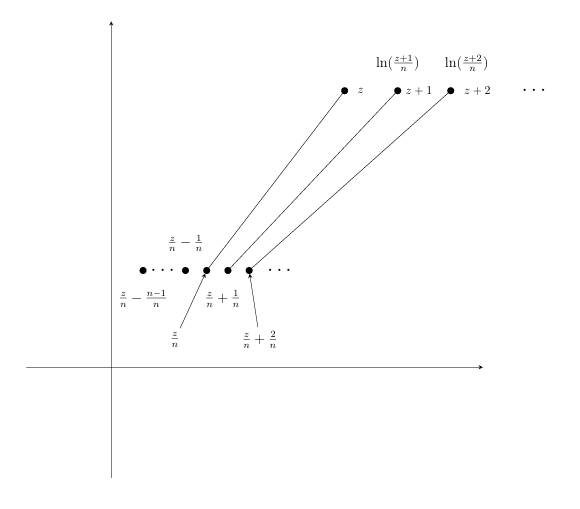}
\\
Since $\ln(n)$ is a constant and $R_{+}[1](z_{0})=-z_{0}-\frac{1}{2}$ it follows from dilation invariance that in equation \ref{GammaMultiplication_3} we have
\begin{equation*}
RHS\,=\,(-z_{0}-\frac{1}{2})\ln n+(z_{0}+\frac{1}{2})\ln n+R_{+}[\ln \tilde{z}](z_{0})\,=\,LHS
\end{equation*}
and the result is proven.\\
\\
\textbf{Comments:}\textbf{(i)} It is worth stopping for a moment to consider the generalised geometric C\'{e}saro definition and treatment of the $\Gamma$-function in this paper, and how it compares with traditional treatments in textbooks and elsewhere. 

One such treatment begins by noticing that $k!=\int_{0}^{\infty}x^{k}e^{-x}\,dx$ and using this to define the gamma function, $\Gamma(z_{0})$, on all of $\mathbb{C}\setminus\mathbb{Z}_{\leq{0}}$ as the Mellin transform of $e^{-x}$, namely $\int_{0}^{\infty}x^{z_{0}-1}e^{-x}\,dx$. The properties of $\Gamma$ such as its functional equation and multiplication formulae are then deduced by applying the theory of contour integration in sophisticated ways - ways which often seem a little opaque in terms of shedding light on why the properties themselves are true. 

Another such method starts by defining $\Gamma(z_{0}+1)$ via the classically convergent Hadamard product formula $\prod_{j=1}^{\infty}(1+\frac{z_{0}}{j}) \cdot \textrm{e}^{-\frac{z_{0}}{j}}$, but again the connection to the original factorial function on positive integers, let alone the more advanced properties of $\Gamma$, is not straightforward. In particular, in all such traditional definitions the proof of the multiplication formulae is rather indirect and obscure, giving little understanding (at least to this author) of \textit{why} these formulae hold.

By contrast, we maintain that the C\'{e}saro definition of the $\Gamma$-function given here (simply as $\prod_{j=1}^{z_{0}}j$) provides as natural as possible an extension of the factorial function from positive integers to all of $\mathbb{C}\setminus\mathbb{Z}_{<0}$. And we maintain that as a consequence it ends up providing a much more direct and insightful interpretation of $\Gamma$ - one in which its properties fall out naturally from the architecture of remainder summation and from the properties of generalised geometric C\'{e}saro convergence in a way which explains transparently why these key properties are true. The above diagrammatic proof of the multiplication formulae as an essentially trivial consequence of the dilation-invariance of C\'{e}saro convergence provides perhaps the best example of this.\\
\\
\textbf{(ii)} To emphasise this point (that the $\Gamma$-function multiplication formulae at heart simply encapsulate the architecture of remainder summation coupled with the dilation-invariance property of C\'{e}saro convergence), note that it is clear from the above proof that an analogous duplication identity will hold for any function defined by remainder summation $R_{+}[f](z)$ (or after exponentiation by a remainder product formula $\prod_{R,+}[g](z)$) as long as $f(n\tilde{z})$ can be related readily to $f(\tilde{z})$. It is trivial for example, by the same reasoning, that
\begin{equation}
\zeta_{H}(z,s)=n^{-s}\sum_{j=1}^{n}\zeta_{H}\left(\frac{z}{n}-\frac{(n-j)}{n},s\right) \quad .
\label{ZetaHDupl}\end{equation}\\
\\  
\textbf{Calculation 2: Proof of our general integral identity for $\zeta_{H}$ :} Earlier in this paper we conjectured the integral identity given in equation  \ref{zeta_Integral_Identity_3}, namely that
\begin{equation*}
\int_{-1}^{0}\,\zeta_{H}\left(z;s\right)\,dz\,=\,0 \quad \textrm{for any} \quad s\in\mathbb{C}\setminus\{1\}.
\end{equation*}
To derive this, suppose initially that $Re(s)<1$. Using right Riemann sums, the integral on the LHS becomes
\begin{equation*}
\int_{-1}^{0}\zeta_{H}(z,s)\,\textrm{d}z=\underset{n\rightarrow\infty}{\lim}\,\frac{1}{n}\,\sum_{j=0}^{n-1}\zeta_{H}(-\frac{j}{n},s)=\underset{n\rightarrow\infty}{\lim}\,\frac{1}{n}\,\sum_{j=0}^{n-1}R_{+}[\tilde{z}^{-s}](-\frac{j}{n})
\end{equation*}
But dilating by n and noting that $(\frac{\tilde{z}}{n})^{-s}=n^{s}\tilde{z}^{-s}$ and that the summation points after dilation all interleave to precisely fill up the positive integers, it follows from C\'{e}saro dilation-invariance that this becomes
\begin{equation*}
\underset{n\rightarrow\infty}{\lim}n^{s-1}R_{+}[\tilde{z}^{-s}](0)=\underset{n\rightarrow\infty}{\lim}n^{s-1}\zeta(s)=0
\end{equation*}
This proves equation (\ref{zeta_Integral_Identity_3}) for $Re(s)<1$. At $s=1$ the integral is not C\'{e}saro convergent but for all other $s$ the result then follows by analytic continuation.\footnote{Note that the reasoning here likewise immediately yields a family of special values of $\zeta_{H}$. For example, for any prime $p$, $\sum\zeta_{H}(-\frac{j}{p},s)=(p^{s}-1)\zeta(s)$ where the sum is over $1<j<p$. For $p,q$ prime, $\sum\zeta_{H}(-\frac{j}{pq},s)=((pq)^{s}-p^{s}-q^{s}+1)\zeta(s)$ where the sum is over $1<j<pq$, j coprime to $pq$, and so on.}

\subsection{Scaling-invariance of C\'{e}saro convergence}

Another natural operator group to consider regarding C\'{e}saro convergence is the scaling group, which is the group of operators $\lbrace S_{r}\rbrace_{r\in\mathbb{R}_{>0}}$ under multiplication ($S_{r}\cdot S_{s}=S_{rs}$) given by $S_{r}[f](z):=f(z^{r})$. Its relationship to C\'{e}saro convergence is more complex than for dilations, but we nonetheless have the same key result:\\
\\
\textbf{Theorem 2:} \textit{Generalised geometric C\'{e}saro convergence is scaling-invariant, i.e. for any} $r\in\mathbb{R}_{>0}$\textit{, if} $g(z):=S_{r}[f](z)=f(z^{r})$ \textit{then} $g$ \textit{has a generalised C\'{e}saro limit along contour $\gamma$ if and only if} $f$ \textit{has a generalised C\'{e}saro limit along contour $\tilde{\gamma}:=\gamma^{r}$ and in this case} $\underset{z\rightarrow\infty}{Clim}\,g(z)\,=\,\underset{z\rightarrow\infty}{Clim}\,f(z)$.\\
\\
\textbf{Proof:} The proof blends the methods of the two proofs of dilation-invariance. To begin with, note that under rescaling $z\mapsto z^{r}$ a C\'{e}saro  eigenfunction/generalised eigenfunction, $z^{\rho} (\ln z)^{m}$, becomes $r^{m}\cdot z^{r\rho} (\ln z)^{m}$, and since $r>0$ (so that $r\rho=0$ if and only if $\rho=0$), it follows that any linear combination of C\'{e}saro eigenfunctions/generalised eigenfunctions with all eigenvalues not equal to $1$ is mapped to a linear combination of such eigenfunctions/generalised eigenfunctions which also has all eigenvalues not equal to $1$. Likewise, on the other hand, any generalised C\'{e}saro eigenfunction with eigenvalue $1$ is preserved by re-scaling. It follows immediately in definition 1 that to prove theorem 3 we need only consider the case where $f(z)$ is C\'{e}saro convergent via a pure power of $P$, namely $P^{n}$ for some $n\in\mathbb{Z}_{>0}$.

But now, to prove this, recall that $P^{-1}=z\frac{d}{dz}+1$ and note that, if we consider $S_{\exp(\epsilon)}$ and again work to order $\epsilon$, then we have that
\begin{equation*}
f(z_{0}^{\exp(\epsilon)})\,=\,f(z_{0}^{(1+\epsilon)})\,=\,f(z_{0}+\epsilon \ln(z_{0})\cdot z_{0})\,=\,f(z_{0})+\epsilon \ln(z_{0}) z_{0}\cdot f^{\prime}(z_{0})+O(\epsilon^{2})
\end{equation*}
so that the generator of re-scalings is
\begin{equation}
H_{S}=z\ln(z)\frac{d}{dz}\,.
\label{Generator_rescalings}\end{equation}\\
Based on this, we have that, for any function $h(z)$,
\begin{eqnarray*}
(P^{-1}\circ H_{S})[h](z) & = & \left(z\frac{d}{dz}+1\right)\left[z\ln(z)\frac{dh}{dz}\right]\\
 & = & z^{2}\ln(z)\frac{d^{2}h}{dz^{2}}+\left(2z\ln(z)+z\right)\frac{dh}{dz}
\end{eqnarray*}
and
\begin{eqnarray*}
(H_{S}\circ P^{-1})[h](z) & = & \left(z\ln(z)\frac{d}{dz}\right)\left[z\frac{dh}{dz}+h(z)\right]\\
 & = & z^{2}\ln(z)\frac{d^{2}h}{dz^{2}} + 2z\ln(z) \frac{dh}{dz} \quad .
\end{eqnarray*}
Thus we have the central commutation relation that
\begin{equation}
P^{-1}\circ H_{S}\,=\,H_{S}\circ P^{-1}+(P^{-1}-1)\,=\,H_{S}\circ P^{-1}+H_{D}
\label{Cesaro_scaling_commutator}\end{equation}
which can alternatively be rewritten as
\begin{equation}
P^{-1}\circ H_{S}\,=\,(H_{S}+1)(P^{-1}-1)+H_{S} \, .
\label{Cesaro_scaling_commutator_2}\end{equation}
This latter form is more convenient since it is then easy to see by iteration that $P^{-1}\circ H_{S}^{2}=(H_{S}+1)^{2}(P^{-1}-1)+H_{S}^{2}$ and $P^{-1}\circ H_{S}^{3}=(H_{S}+1)^{3}(P^{-1}-1)+H_{S}^{3}$ and in general
\begin{equation}
P^{-1}\circ H_{S}^{n}\,=\,(H_{S}+1)^{n}(P^{-1}-1)+H_{S}^{n}
\label{Cesaro_scaling_commutator_3}\end{equation}
for any $n\in\mathbb{Z}_{>0}$. But then it follows immediately from equation \ref{Cesaro_scaling_commutator_3} and from the fact that $S_{r}=e^{rH_{s}}$ that we have, in general,
\begin{eqnarray*}
P^{-1}\circ S_{r} & = & \sum_{n=0}^{\infty} \frac{r^{n}}{n!}P^{-1}\circ H_{S}^{n} \,=\, \sum_{n=0}^{\infty} \frac{r^{n}}{n!}\lbrace(H_{S}+1)^{n}(P^{-1}-1)+H_{S}^{n}\rbrace\\
 & = & e^{r(H_{s}+1)}(P^{-1}-1)+e^{rH_{s}} \,=\, e^{r}S_{r}(P^{-1}-1)+S_{r}\, .
\end{eqnarray*}
Pre- and post-multiplying by $P$ and re-arranging we thus obtain the following quasi-commutation relation:
\begin{equation}
P \circ S_{r}\,=\,\tilde{q}(P) \circ S_{r}\circ P \quad where \quad \tilde{q}(P)=\left(\frac{e^{r}-1}{e^{r}}\right)P+\frac{1}{e^{r}} \, .
\label{Cesaro_scaling_commutator_4}\end{equation}
Since $\tilde{q}(P)$ is obviously regular (i.e. $\tilde{q}(1)=1$), it follows at once that if $P[f]\rightarrow L$ classically then $g(z):=S_{r}[f](z)$ also satisfies that $P[g]\rightarrow L$ classically as $z\rightarrow\infty$. And the converse likewise holds since we can equally view $f$ as $S_{\frac{1}{r}}[g]$. As for the case of higher powers of $P$, i.e. where we need to take $P^{n}[f]$ for some $n\in\mathbb{Z}_{>1}$ in order to get a classically convergent function, note that iterating equation \ref{Cesaro_scaling_commutator_4} leads immediately to the generalised commutation relation that, for any such $n$
\begin{equation}
P^{n} \circ S_{r}\,=\,\tilde{q}(P)^{n} \circ S_{r} \circ P^{n} \quad .
\label{Cesaro_scaling_commutator_5}\end{equation}
It thus follows at once by identical reasoning that $P^{n}[f]\rightarrow L$ classically as $z\rightarrow\infty$ if and only if $g(z):=S_{r}[f](z)$ also satisfies that $P^{n}[g]\rightarrow L$ classically. This completes the proof of theorem 2.\\
\\
\textbf{Comment:} Theorem 2 seems an intriguing result. At this stage, however, we have not found the same sort of direct applications of it, either to existing or new problems, that we previously illustrated regarding the dilation-invariance of generalised C\'{e}saro convergence. As such, we do not explore scaling-invariance any further here, other than to encourage others to investigate such possible applications; and this in turn concludes our working in this paper.

\section{Appendices}

We include as appendices two brief additional calculations. The first, as promised in the final comment in section 2.3, is the derivation of the functional equation for $\zeta(s)$ from consideration of bi-directional remainder summation in $z_{0}$ for $\zeta_{H}(z_{0};s)$. The second is a quick observation facilitating the rapid numerical calculation of $\Gamma(z_{0})$ using the C\'{e}saro definition we have introduced in this paper.

\subsection{Derivation of the functional equation for $\zeta(s)$ from bi-directional summation for $\zeta_{H}(z_{0};s)$}

Recall comments (i) and (ii) following the derivation of the functional equation for $\Gamma$ in section 2.3. If we seek a functional equation for $\zeta_{H}(z_{0};s)$ - or rather for $\tilde{\zeta}_{H}(z_{0};s):=R_{+,0}[\tilde{z}^{-s}](z_{0};s)=\zeta_{H}(z_{0};s)+z_{0}^{-s}$ - then we are naturally led to consider the bi-directional sum $R_{+,0,-}[\tilde{z}^{-s}](z_{0};s)$ and its Fourier series representation as a periodic function of $z_{0}$ with period $1$.

If we write this Fourier series representation as
\begin{equation*}
R_{+,0,-}[\tilde{z}^{-s}](z_{0})=\sum_{n=-\infty}^{\infty}a_{n}e^{2\pi inz_{0}}
\end{equation*}
then we have that
\begin{equation*}
a_{n}=\int_{0}^{1}R_{+,0,-}[\tilde{z}^{-s}](z_{0};s)e^{-2\pi inz_{0}}\,dz_{0}=\int_{-\infty}^{\infty}z^{-s}e^{-2\pi inz}\,dz = \mathcal{F}[\tilde{z}^{-s}](n)
\end{equation*}
Now, recall that $\mathcal{F}[1](\xi)=\delta(\xi)$ and that $\mathcal{F}[zf(z)](\xi)=\frac{1}{-2\pi i}\frac{d}{d\xi}\hat{f}(\xi)$. It follows that $\mathcal{F}[\frac{1}{\tilde{z}}](\xi)=-2\pi i u(\xi)$, and $\mathcal{F}[\frac{1}{\tilde{z}^{2}}](\xi)=(-2\pi i)^{2} \xi_{+}$ (where $\xi_{+}$ equals $\xi$ if $\xi \geq{0}$ and equals $0$ if $\xi<0$ in the usual way), and $\mathcal{F}[\frac{1}{\tilde{z}^{3}}](\xi)=(-2\pi i)^{3} \frac{\xi_{+}^{2}}{2!}$; and in general $\mathcal{F}[\tilde{z}^{-s}](\xi)=(-2\pi i)^{s} \frac{\xi_{+}^{s-1}}{\Gamma(s)}=e^{-i\pi s/2}2^{s}\pi^{s} \frac{\xi_{+}^{s-1}}{\Gamma(s)}$. Thus, taking $Re(s)>1$ initially in order to avoid technicalities for $n=0$, we have that
\begin{equation*}
a_{n}= \begin{cases}
e^{-i\pi s/2}\frac{2^{s}\pi^{s}}{\Gamma(s)}n^{s-1} \quad, & n\geq{0}\\
0 \quad, & n<0\end{cases}
\end{equation*}
and it follows that
\begin{equation}
R_{+,0,-}[\tilde{z}^{-s}](z_{0})=e^{-i\pi s/2}\frac{2^{s}\pi^{s}}{\Gamma(s)}\sum_{n=1}^{\infty}n^{s-1}e^{2\pi inz_{0}} \quad .
\label{Bidirectional_zeta_1}\end{equation}
This is the general expression from which a functional equation relating the values of $\tilde{\zeta}_{H}(z_{0};s)$ and $\tilde{\zeta}_{H}(1-z_{0};s)$ could then be derived in the manner previously discussed. Our focus here, however, is just on deriving the functional equation in $s$ for the Riemann zeta function $\zeta(s)$. For this we simply take equation \ref{Bidirectional_zeta_1} and consider the particular case where $z_{0}=\frac{1}{2}$. Then we get that
\begin{eqnarray*}
LHS & = & \left\{ \ldots+\left(\frac{-3}{2}\right)^{-s}+\left(\frac{-1}{2}\right)^{-s}+\left(\frac{1}{2}\right)^{-s}+\left(\frac{3}{2}\right)^{-s}+\ldots \right\}\\
 & = & \left(1+e^{-i\pi s}\right)\left\{ \left(\frac{1}{2}\right)^{-s}+\left(\frac{3}{2}\right)^{-s}+\ldots \right\}\\
 & = & \left(1+e^{-i\pi s}\right)2^{s}\left(1-2^{-s}\right)\zeta(s) \,=\, \left(1+e^{-i\pi s}\right)\left(2^{s}-1\right)\zeta(s)
\end{eqnarray*}
and
\begin{eqnarray*}
RHS & = & e^{-i\pi s/2}\frac{2^{s}\pi^{s}}{\Gamma(s)}\sum_{n=1}^{\infty}(-1)^{n}n^{s-1}\\
 & = & e^{-i\pi s/2}\frac{2^{s}\pi^{s}}{\Gamma(s)}(2^{s}-1)\zeta(1-s) \, ,
\end{eqnarray*}
where we have used dilation-invariance in the course of the C\'{e}saro treatment of the divergent series involved. It follows that
\begin{equation*}
\left(1+e^{-i\pi s}\right)\zeta(s)= e^{-i\pi s/2}\frac{2^{s}\pi^{s}}{\Gamma(s)}\zeta(1-s)
\end{equation*}
and after elementary re-arrangement this yields the desired functional equation for $\zeta$, namely that
\begin{equation}
\zeta(1-s)=2^{1-s}\pi^{-s}\cos(\frac{\pi s}{2})\Gamma(s)\zeta(s) \, .
\label{Zeta_functional_eqn}\end{equation}
\textbf{Comment:} Note that, while the general approach we previously outlined explains why we introduce the bi-directional sum in this working (and explains why the functional equation for $\tilde{\zeta}_{H}$ would relate function values at $z_{0}$ and $1-z_{0}$ had we derived it in general), the reason that the functional equation in $s$ for $\zeta$ relates $\zeta(s)$ and $\zeta(1-s)$ is for a different reason. It is because of the Fourier-transform relationship that maps $z^{-s}$ to $\xi^{s-1}$, not because of any hidden bi-directional summation in $s$ (as far as we can tell!). Nonetheless, this calculation once again confirms the utility of generalising from $\zeta(s)$ to $\zeta_{H}(z_{0};s)$ and then considering bi-directional summation in $z_{0}$ as a mechanism for obtaining a periodic function - which allows us in turn to invoke Fourier theory and hence obtain the functional equation for zeta by particularising to $z_{0}=\frac{1}{2}$.\footnote{When dealing with problems involving remainder summation and the generalised C\'{e}saro framework, it always seems to be a useful thing to also consider bi-directional summation!}

\subsection{Efficient numerical calculation of $\Gamma(z_{0}+1)$ using its C\'{e}saro definition}

Recall, from equations \ref{Remprodzb} and \ref{ln_Gamma_tilde_formula} in the first calculation of section 2.3, that the log-Gamma function is given by
\begin{equation*}
\ln(\Gamma(z_{0}+1))= \frac{1}{2}\ln(2\pi)-C_{z_{0}}
\end{equation*}
where
\begin{equation*}
C_{z_{0}}=\underset{k\rightarrow\infty}{\lim}\left\{ \sum_{j=1}^{k}\ln(z_{0}+j)-(z_{0}+k+\frac{1}{2})\ln(z_{0}+k)+(z_{0}+k)\right\}
\end{equation*}
Theoretically this allows numerical calculation of $\ln(\Gamma(z_{0}+1))$ (and hence immediately of $\Gamma(z_{0}+1)$) for arbitrary $z_{0}\in\mathbb{C}\setminus\mathbb{Z}_{<0}$ by taking enough terms in the limiting sum defining $C_{z_{0}}$ to achieve the desired accuracy - but the convergence is slow! 

In deriving this formula for $C_{z_{0}}$ in section 2.3, however, we only considered sufficiently many terms in the Euler-McLaurin sum formula for $\sum_{j=1}^{k}\ln(z_{0}+j)$ to achieve convergence, without concerning ourselves at all with speed or numerical efficiency. If we simply include additional terms in this Euler-McLaurin expansion then, subject to all the usual caveats regarding asymptotic series, we can obtain alternative, faster expressions for calculating $C_{z_{0}}$ and these permit extremely rapid numerical calculation of $\ln(\Gamma(z_{0}+1))$ for arbitrary $z_{0}$ with extremely simple code. For example, if we keep an additional three terms in the Euler-Mclaurin expansion we obtain
\begin{equation}
C_{z_{0}}=\underset{k\rightarrow\infty}{\lim}\left\{\begin{array}{cc}\sum_{j=1}^{k}\ln(z_{0}+j)-(z_{0}+k+\frac{1}{2})\ln(z_{0}+k)+(z_{0}+k)\\
\\
-\frac{1}{12}\frac{1}{z_{0}+k}+\frac{1}{360}\frac{1}{(z_{0}+k)^{3}}-\frac{1}{1260}\frac{1}{(z_{0}+k)^{5}}\end{array}\right\}
\label{Numerical_ln_Gamma_sample_1}\end{equation}
and similar formulae can be trivially obtained by taking fewer or more terms in this expansion as desired.

\section{Acknowledgements}

We thank Professor John Keats and Professor Robert Frost (pers. comm.) for many insightful comments, Professor T. Abby for his help in preparing this paper, and Christiana Stone for her help in creating the figures.


\begin{thebibliography}{4}
\bibitem{key-1}R. Stone, \textit{Introduction to generalised C\'{e}saro convergence I}, 2026

\bibitem{key-2}Lars V. Ahlfors, \textit{Complex Analysis: An Introduction to the Theory of Analytic Functions of One Complex Variable}, Third Edition, McGraw-Hill, 1979

\bibitem{key-3}H.M. Edwards, \textit{Riemann's Zeta Function}, Academic Press, 1974 
\end{thebibliography}
\end{document}